\documentclass[a4paper,12pt]{article}
\usepackage{epsf,amsfonts,url}
\def\g{\gamma}
\def\R{\mathbb{R}}
\def\Z{\mathbb{Z}}
\def\CP{\mathbb{C}P}

\def\thep{.}
\newcommand{\beqn}{\begin{equation}}
\newcommand{\eeqn}{\end{equation}}
\newcommand{\beqna}{\begin{eqnarray}}
\newcommand{\beqnao}{\begin{eqnarray*}}
\newcommand{\eeqna}{\end{eqnarray}}
\newcommand{\eeqnao}{\end{eqnarray*}}
\newcommand{\ba}{\begin{array}}
\newcommand{\ea}{\end{array}}
\newtheorem{Theorem}{Theorem}[section]
\newtheorem{lemma}[Theorem]{Lemma}

\newtheorem{definition}[Theorem]{Definition}

\newcommand{\one}{\relax\ifmmode{\rm 1\>\!\!\!I}\else{$\rm 1\!I$}\fi}

\begin{document}
\def\q{Q}

\title{Tri-hamiltonian Toda lattice and a canonical bracket for closed discrete curves}
\author{ Nadja Kutz\footnote{Supported by the Servicezentrum Humboldt
   University of Berlin.
    email: \protect\url{nadja@math.tu-berlin.de}}
  }
\maketitle
\label{sec:dc}
\begin{abstract}
Flows on (or variations of) discrete curves in  $\R^2$  give rise to flows on
a subalgebra of functions on that curve. 
For a special choice of flows and  a certain subalgebra this is
described by the Toda lattice hierachy \cite{HK02}. 
In the paper it is shown that the canonical 
symplectic structure on
$\R^{2N}$, which can be interpreted as the phase space  
 of closed discrete curves in $\R^2$ with length $N,$ 
induces Poisson commutation
relations on the above mentioned subalgebra which yield the 
tri-hamiltonian poisson structure
of the Toda lattice hierachy.
\end{abstract}
\section{Introduction}The Toda lattice hierarchy is a set of
equations, including the Toda equation:
$$
\ddot q_k = e^{q_{k+1}-q_k}-e^{q_{k}-q_{k-1}}.
$$
The Toda equation is sometimes also called first flow equation of the
Toda lattice hierarchy, it was discovered by Toda in 1967 \cite{TO89}. 
A good overview about the  literature about the Toda lattice
can be found in \cite{FT86,SU02}.
It is a wellknown fact that the Toda lattice hierachy has a socalled 
trihamiltonian structure (for an overview about trihamiltonian
structures please
see \cite{SU02} and the therein cited references). In particular 
this means that the Toda lattice as a dynamical system admits
three different poisson structures. In the following article it
will be shown explicitly how these three different structures can be 
globally obtained
from the the canonical poisson structure on $\R^{2N}$. We will confine ourselves
to the study of the periodic Toda lattice hierachy. To do so
we define a map from canonical 
coordinates on  $\R^{2N}$ to
the $2N$ (periodic) Flaschka-Manakov \cite{MAN74,Flasch74a,Flasch74b} 
variables that depends on a spectral parameter. The canonical poisson
structure on $\R^{2N}$ induces then a Poisson structure for the
 periodic Flaschka-Manakov variables, which  is the tri-hamiltonian
Poisson structure of the Toda lattice .

The paper is organized as follows: We will briefly recall the
connection between discrete curves in $\R^2$ and the Toda lattice
hierachy, for a more thorough investigation of that connection 
please see \cite{HK02}. In particular this brief part should
serve as a motivation for how the spectral parameter dependent
map is derived. The phase space of all closed (real) discrete curves
of length $N$ is  $\R^{2N}$, hence the canonical poisson structure 
on $\R^{2N} $ is a Poisson
structure on the phase space of all closed discrete curves.
In the second part we will then state that this canonical Poisson structure 
is in fact a Poisson structure which gives the three brackets
of the (periodic) Toda lattice.
\section{Discrete curves and the Toda lattice}
\begin{definition}
A discrete curve in $\R^2$ is  a map 
\beqna \g:\Z&\to&\R^2\nonumber\\
k&\mapsto& \g_k = \left(\ba{c} x_k \\ y_k \ea\right)
\label{def:dc:gamdef}\eeqna
\end{definition}

\noindent Define:
\beqna
g_k &=& \det(\g_k,\g_{k+1}) = x_k y_{k+1}-y_k x_{k+1}\\
u_k&=&\det(\g_{k-1},\g_{k+1})=x_{k-1}y_{k+1}-y_{k-1} x_{k+1} \label{def:dc:gu}
\eeqna
We will from now on consider the generic case $g_k\neq0$ for
all $k$.
The following lemma can be straightforwardly obtained by using
the above definitions (\ref{def:dc:gu}):
\begin{lemma}
\beqn\g_{k+1}=\frac{1}{g_{k-1}}(u_k \g_k - g_k \g_{k-1}).
\label{eq:dc:GamRec}
\eeqn
\label{lem:dc:Gam1}\end{lemma}
If the  variables $u_k$ and $g_k$ and initial points
$\g_0$ and $\g_1$ are given, then lemma~\ref{lem:dc:Gam1} is
a recursive definition of a discrete curve.

Note that $\det(\g_k, \frac{\g_{k+1}-\g_{k-1}}{g_k + g_{k-1}})= 1$. This
means that  
$\g_k$ and  $\frac{\g_{k+1}-\g_{k-1}}{g_k + g_{k-1}}$
 are linearly independent. 
Hence an arbitrary flow on (or variation of) $\g$  can be written in the 
following way:
\begin{equation}
       \frac{d}{dt}\g_k= \dot\g_k = \alpha_k \g_k + \frac{\beta_k}{u_k}(\g_{k+1}-\g_{k-1})
       \hspace{4mm}\alpha_k,\beta_k \in \R \thep\label{eq:dc:generalFlow}
\end{equation}
The variables $\alpha_k,\beta_k \in \R $ are arbitrary.
Equations (\ref{eq:dc:GamRec}) and (\ref{eq:dc:generalFlow}) can also be reformulated as a zero-curvature condition.
Define:
\beqn
F_k
=\left(\ba{l}\g^T_{k}\\\g^T_{k-1}\ea\right)=
\left(\ba{ll}x_k&y_k\\x_{k-1}&y_{k-1}\ea \right).
\label{def:dc:Fdef}\eeqn
\begin{lemma}
Let $\alpha_k, \beta_k \in \R$ be arbitrary and $g_k$,$u_k$ be as
defined in (\ref{def:dc:gu}). Then 
$$F_{k+1}=L_k F_k \hspace{5mm} \dot F_k = V_k F_k$$
with
\beqn
L_k = \left(
  \ba{cc}\frac{1}{g_{k-1}}u_k&-\frac{g_k}{g_{k-1}}\\1&0\ea\right)
\hspace{5mm}
V_k = \left(\ba{cc}
\alpha_{k}+\frac{1}{g_{k-1}}\beta_{k}&
-(1+\frac{g_{k}}{g_{k-1}})\frac{\beta_{k}}{u_{k}}\\
(1+\frac{g_{k-2}}{g_{k-1}})\frac{\beta_{k-1}}{u_{k-1}}&
\alpha_{k-1}-\frac{1}{g_{k-1}}\beta_{k-1}
\ea \right).
\label{eq:dc:vk}\eeqn
The compatibility equation
\beqn
\dot L_k = V_{k+1}L_k - L_k V_k \label{eq:dc:comp}
\eeqn
is satisfied for all $\alpha_k, \beta_k \in \R$. 
\label{prop:dc:LV}
\end{lemma}
The compatibility equation (\ref{eq:dc:comp}) is also called
zero curvature equation or condition. We call $F_k$ a discrete frame.

\begin{lemma}
 A flow on the discrete curve $\g$ given by (\ref{eq:dc:generalFlow}) 
 generates the following flow on the variables $g_k$:
\begin{equation}
       \dot g_k = g_k (\alpha_{k+1}+\alpha_k) +\beta_{k+1}-\beta_k\thep
        \label{eq:dc:dotg}
\end{equation}
\label{lem:dc:gu}\end{lemma}

\begin{lemma}
A flow on the discrete curve  $\g$ given by (\ref{eq:dc:generalFlow}) 
 generates the following flow on the variables $u_k$:
\beqn\begin{array}{rcl}
\dot u_k &=& u_k(\alpha_{k-1}+ \alpha_{k+1})\\&+&
 \beta_{k-1}\frac{g_{k}}{u_{k-1} g_{k-1}}(g_{k-2}+g_{k-1}) -
\beta_{k+1}\frac{g_{k-1}}{u_{k+1} g_{k}}(g_{k}+g_{k+1})\\
&+&u_k (\frac{1}{g_{k}}\beta_{k+1} - \frac{1}{g_{k-1}}\beta_{k-1}) 
\label{eq:dc:FlowOnU}
\end{array}
\eeqn
\end{lemma}
Define
\beqna
a_k &:=& g_k^{-2}\nonumber\\
b_k &:=& \frac{u_k}{g_{k-1}g_{k}} - \lambda,
\label{def:dc:abdef}\eeqna
where $\lambda$ is an arbitrary (but in particular time independent) parameter.
Clearly the flows on the variables $g_k$ and $u_k$ given in
(\ref{eq:dc:dotg}) and (\ref{eq:dc:FlowOnU}) define flows on
the variables $a_k$, $b_k$ via definitions (\ref{def:dc:abdef}).
\begin{Theorem}
Denote
$$
V_k := \left( \ba{cc}v_k^{11}&v_k^{12}\\v_k^{21}&v_k^{22}\ea\right).
$$
Define
\beqn
\alpha_k := v_k^{11}+\frac{v_k^{12}u_k}{g_{k-1}+g_k}
\hspace{5mm}
\beta_k := -\frac{v_k^{12} g_{k-1}u_k}{g_{k-1}+g_k}.
\label{def:dc:alphabet}
\eeqn
By (\ref{eq:dc:generalFlow}) and with definition
(\ref{def:dc:alphabet}),
 $\alpha_k$ and 
$\beta_k$ define a  flow on discrete curves in
$\R^2$ depending on the choice of the $V_k$. This induces via
the definition (\ref{def:dc:abdef}) a flow on the 
variables $a_k,b_k$. 

If the $\{V_k\}_{k\in\{1...N\}}$ ($V_k$ periodic) are the 
two dimensional Lax matrices defining 
the n'th Toda flow (see e.g.\cite{FT86,SU02}) with $\lambda$ being the
corresponding spectral parameter, then 
$\{a_k,b_k\}_{k\in\{1...N\}}$ ($a_k,b_k$ periodic) are the
Flaschka-Manakov variables obeying the n'th Toda flow.
\label{theor:dc:Toda}\end{Theorem}
For a proof of theorem (\ref{theor:dc:Toda}) and lemmas 
(\ref{eq:dc:FlowOnU}) and (\ref{eq:dc:dotg}) please see \cite{HK02}
Note that the Toda flows (as flows on the 
Flaschka-Manakov variables $a_k,b_k$) do not depend
on the choice of $\lambda$, whereas the corresponding flows for
the determinants $g_k$ and $u_k$ of course depend on the choice of $\lambda$.

\section{How to derive Toda Poisson  brackets  from the  canonical Poisson bracket}
The following theorem is easy to state and can be verified straightforwardly.
Nevertheless the assertion itself is  not really suggestive apriori.  
And indeed it was not found by good guessing and then verifying, but 
by starting with the
brackets of the Toda model (while assuming that they can be derived
from the canonical coordinates via the below map) and the help  of a 
computer algebra system.
\begin{Theorem}
Let $N>3$. Let $\{x_i,y_i | i \in \{0...N-1\}\}$ be canonical 
coordinates on $\R^{2N}$. Define the following symplectic structure on  
$\R^{2N}$:
$$
\Omega := 2 \sum_{i=0}^{N} dx_i \wedge dy_i
$$
which leads to the ultralocal Poisson relations:
$$
\{x_i,y_i\}= \frac{1}{2}, \hspace{1cm} \mbox{and zero else.}
$$
Let $x,y$ describe a periodic phase space, i.e:
\beqna
x_{k+N}&:=&x_k\\
y_{k+N}&:=&y_k
\eeqna
Define:
\beqna
g_k &:=& x_k y_{k+1} - y_k x_{k+1}\\
u_k &:=& x_{k-1} y_{k+1} - y_{k-1} x_{k+1}\\
a_k &:=& g_{k}^{(-2)}\\
b_k &:=& \frac{u_{k}}{g_{k - 1}g_{k}} - \lambda,
\eeqna
where $\lambda \in \R$ arbitrary.

Then the Poisson relations for the variables $a_i$ and $b_i$, which are
given via the ultralocal relations for $x_i$ and $y_i$ read as:
\beqna
\{a_k,a_{k+1}\}&=&-2 a_k a_{k+1} b_{k+1}-2 a_k a_{k+1}\lambda\label{def:dc:CR1}\\
\{b_k,b_{k+1}\}&=&-a_k (b_k + b_{k+1}) - 2 a_k \lambda\label{def:dc:CR2}\\
\{b_k,a_{k-2}\}&=&a_{k-2}a_{k-1}\\
\{b_k,a_{k-1}\}&=&a_{k-1}(b_{k}^2+a_{k-1})+2 b_k a_{k-1} \lambda + a_{k-1} \lambda^2\label{def:dc:CR3}\\
\{b_k,a_{k}\}&=&-a_k(b_k^2+a_k)-2 b_k a_k \lambda - a_k \lambda^2 \label{def:dc:CR4}\\
\{b_k,a_{k+1}\}&=&-a_k a_{k+1}\label{def:dc:CR5}
\eeqna
and zero for all the remaining commutators. For each power 
of $\lambda$ this gives one of the three Poisson 
brackets of the Toda lattice hierachy (as e.g. given in \cite{SU02}).
\end{Theorem}
Note that the brackets (\ref{def:dc:CR1})-(\ref{def:dc:CR5}) make only sense
if $N>3$. Nevertheless the canonical Poisson structure 
defines also a Poisson structure on the Flaschka-Manakov variables
for N=1,2,3.
\section{Conclusions}
The extremly nice Poisson bracket for the coordinates of closed discrete
curves, which lead to the brackets of the periodic Toda lattice suggests
that the study of discrete curves in the context of the Toda lattice
is not an exotic geometrical excursion but rather a very natural
choice.
 
The map from the canonical coordinates $x_k$,$y_k$ to the Flaschka-Manakov
variables $a_k$, $b_k$ is not one-to-one. This becomes also clear by studying
the phase space of quasiperiodic discrete curves, i.e. curves for which
$\g_{k+N}= {\cal M}\g_{k},$ where ${\cal M}\in GL(2,\R)$. If 
${\cal M}\in SL(2,\R)$ then the corresponding Flaschka-Manakov
variables are also periodic, but are  not necessarily restricted to the leaves
given by the Casimir functions. Nevertheless we expect that  
Poisson brackets for that general case will be rather nontrivially.
We are investigating that phase space right now.  

The canonical coordinates are matrix entries of the discrete frame $F_k$ given in (\ref{def:dc:Fdef}). In that sense the canonical bracket 
can be seen as a bracket for the frame. The frame is  the
(partially discrete) integral of the zero curvature condition in 
(\ref{eq:dc:comp})
subject to a periodicity condition (i.e. being a closed curve). It is
an interesting question wether other integrable systems are admitting
such a simple bracket for their corresponding frames. 

So it may be also important to find out, wether there is a connection 
to the work of Gelfand
and Zakharevich \cite{GZ},
where its proven that there exists a local isomorphism between
the bihamiltonian bracket of the periodic Toda lattice and a product of
two canonical brackets (at generic points).

It would also be interesting to study wether the above could be
extended to generalized Toda systems \cite{Ko79}.

\section{Acknowledgements}
The author would like to thank Tim Hoffmann for 
helpful discussions.


\end{document}